\documentclass[12pt]{article}

\usepackage{amsmath,amssymb}
\usepackage{color}

\newtheorem{theorem}{Theorem}[section]

\newtheorem{corollary}[theorem]{Corollary}

\newtheorem{observation}[theorem]{Observation}

\newcommand{\proof}{\noindent{\bf Proof.\ }}
\newcommand{\qed}{\hfill $\square$ \bigskip}
\def\cp{\,\square\,}
\DeclareMathOperator {\gp} {gp}

\DeclareMathOperator {\diam} {diam}

\begin{document}

\title{Characterization of general position sets and its applications to cographs and bipartite graphs}

\author{
    Bijo S. Anand $^{a}$
    \and
	Ullas Chandran S. V. $^{b}$
	\and
	Manoj Changat $^{c}$
	\and
	Sandi Klav\v zar $^{d,e,f}$
	\and
	Elias John Thomas $^{g}$ 
}

\date{\today}

\maketitle
\begin{center}
$^a$ Department of Mathematics, Sree Narayana College, Punalur-691305, Kerala, India;\quad {\tt bijos\_anand@yahoo.com} 
\medskip

$^b$ Department of Mathematics, Mahatma Gandhi College, Kesavadasapuram,  Thiruvananthapuram-695004, Kerala, India;\quad 
{\tt svuc.math@gmail.com} 
\medskip

$^c$ Department of Futures Studies, University of Kerala Thiruvananthapuram-695034, Kerala, India;\quad  
{\tt mchangat@gmail.com} 
\medskip

$^d$ Faculty of Mathematics and Physics, University of Ljubljana, Slovenia\\
{\tt sandi.klavzar@fmf.uni-lj.si}
\medskip

$^e$ Faculty of Natural Sciences and Mathematics, University of Maribor, Slovenia
\medskip

$^f$ Institute of Mathematics, Physics and Mechanics, Ljubljana, Slovenia
\medskip

$^g$ Department of Mathematics, Mar Ivanios College, Thiruvananthapuram-695015, Kerala, India;\quad 
{\tt eliasjohnkalarickal@gmail.com}
		
\end{center}

\begin{abstract}
A vertex subset $S$ of a graph $G$ is a general position set of $G$ if no vertex of $S$ lies on a geodesic between two other vertices of $S$.  The cardinality of a largest general position set of $G$ is the general position number ${\rm gp}(G)$ of $G$. It is proved that $S\subseteq V(G)$ is in general position if and only if the components of $G[S]$ are complete subgraphs, the vertices of which form an in-transitive, distance-constant partition of $S$. If ${\rm diam}(G) = 2$, then ${\rm gp}(G)$ is the maximum of $\omega(G)$ and the maximum order of an induced complete multipartite subgraph of the complement of $G$. As a consequence, ${\rm gp}(G)$ of a cograph $G$ can be determined in polynomial time. If $G$ is bipartite, then ${\rm gp}(G) \leq \alpha(G)$ with equality if ${\rm diam}(G) \in \{2,3\}$. A formula for the general position number of the complement of an arbitrary bipartite graph is deduced and simplified for the complements of trees, of grids, and of hypercubes. 
\end{abstract}


\medskip\noindent
{\bf AMS Subj.\ Class:} 05C12; 05C69; 68Q25.

\section{Introduction}
\label{sec:intro}

Motivated by the century old Dudeney's no-three-in-line problem~\cite{dudeney-1917}
(see~\cite{ku-2018, misiak-2016, PoWo07} for recent developments on it) and by the general position subset selection problem~\cite{froese-2017, payne-2013} from discrete geometry, the natural related problem was introduced to graph theory in~\cite{manuel-2018a} as follows. Let $G=(V(G), E(G))$ be a graph. Then we wish to find a largest set of vertices $S\subseteq V(G)$, called a {\em gp-set} of $G$, such that no vertex of $S$ lies on a geodesic (in $G$) between two other vertices of $S$. The {\em general position number} ({\em $\gp$-number} for short), $\gp(G)$, of $G$ is the number of vertices in a gp-set of $G$. 

As it happens, the same concept has already been studied two years earlier in~\cite{ullas-2016} under the name {\em geodetic irredundant sets}. The concept was formally defined in a different, more technical language, see the preliminaries below.  In~\cite{ullas-2016} graphs $G$ with $\gp(G)\in \{2, n(G)-1, n(G)\}$ were characterized and several additional results about the general position number were deduced. The term general position problem was coined in~\cite{manuel-2018a}, where different general upper and lower bounds on the gp-number are proved. In the same paper it is demonstrated that in a block graph the set of simplicial vertices forms a gp-set and that the problem is NP-complete in the class of all graphs. In the subsequent paper~\cite{manuel-2018b} the gp-number is determined for a large class of subgraphs of the infinite grid graph, for the infinite diagonal grid, and for Bene\v{s} networks.

In this paper we continue the investigation of general position sets in graphs. In the following section definitions and preliminary observations are listed. In Section~\ref{sec:characterization} we prove a characterization of general position sets and demonstrate that some earlier results follow directly from the characterization. In the subsequent section we consider graphs of diameter $2$. We prove that if $G$ is such a graph, then $\gp(G)$ is the maximum of the clique number of $G$ and the maximum order of an induced complete multipartite subgraph of the complement of $G$. In the case of cographs the latter invariant can be replaced by the independence number. As a consequence, $\gp(G)$ of a cograph $G$ is polynomial. Moreover, we determine a formula for $\gp(G)$ for graphs with at least one universal vertex. In Section~\ref{sec:bipartite} we consider bipartite graphs and their complements. If $G$ is bipartite, then $\gp(G) \leq \alpha(G)$ with equality if $\diam(G) \in \{2,3\}$.  We prove a formula for the general position number of the complement of a bipartite graph and simplify it for the complements of trees, of grids, and of hypercubes. In particular, $\gp(\overline{T}) = \max \{\alpha(T), \Delta(T) +1\}$ for a tree $T$.

\section{Preliminaries}
\label{sec:preliminary}

Let $G$ be a connected graph and $u,v\in V(G)$. The {\it distance} $d_G(u,v)$ between $u$ and $v$ is the minimum number of edges on a $u, v$-path.  The maximum distance between all pairs of vertices of $G$ is the {\em diameter} $\diam(G)$ of $G$. A $u,v$-path of length $d_G(u,v)$ is called an $u,v$-{\it geodesic}. The {\em interval} $I_G[u,v]$ between vertices $u$ and $v$ of a graph $G$ is the set of vertices $x$ such that there exists a $u,v$-geodesic which contains $x$. For $S\subseteq V(G)$ we set $I_G[S]=\bigcup_{_{u,v\in S}}I_G[u,v]$. To simplify the writing, we may omit the index $G$ in the above notation provided that $G$ is clear from the context.

A set of vertices $S\subseteq V(G)$ is a {\em general position set} of $G$ if no three vertices of $S$ lie on a common shortest path. A gp-set is thus a largest general position set. Call a vertex $v\in T\subseteq V(G)$ to be an {\em interior vertex} of $T$, if $v \in I [T - \{v\}]$. Now, $T$ is a general position set if and only if $T$ contains no interior vertices. In this way general position sets were introduced in~\cite{ullas-2016} under the name geodetic irredundant sets. 

The set $S$ is \emph{convex} in $G$ if $I[S]=S$. The \emph{convex hull} $H(S)$ of $S$ is the smallest convex set that contains $S$, and $S$ is a \emph{hull set} of $G$ if $H(S)$ is the whole vertex set of $G$. A  smallest hull set is a \emph{minimum hull set} of $G$, its cardinality is the \emph{hull number} $h(G)$ of $G$. A hull set $S$ in a graph $G$ is a \emph{minimal hull set} if no proper subset of $S$ is a hull set of $G$. The number of vertices in a largest minimal hull set of $G$ is its \emph{upper hull number} $h^+(G)$. It is clear that in any graph $G$, every minimum hull set is also minimal, therefore $h(G)\leq h^+(G)$. The following fact is obvious. 

\begin{observation}
\label{obs1} 
Let $G$ be a connected graph, $S$ a minimal hull set, $x, y\in S$, and $P$ a $x,y$-geodesic. If $z\in P$, where $z\ne x,y$, then $z\notin S$.
\end{observation} 

It follows from this observation that every minimal hull set is a  general position set. Consequently, 
$$2 \leq h(G) \leq h^+(G) \leq \gp(G) \leq n(G)\,,$$ 
where $n(G) = |V(G)|$. The related difference can be arbitrary large though. For instance, if $n\ge 2$, then $h^+(K_{n,n}) = 2$ and $\gp(K_{n,n}) = n$. 

With respect to convexities we mention the following parallel concept to the general position number, where in the definition of the interior vertex we replace ``I'' with ``H''. More precisely, the {\em rank} of a graph $G$ is the cardinality of a largest set $S$ such that $v\notin H(S - \{v\})$ for every $v\in S$, see~\cite{kante-2017}. Actually, the graph rank can be studied for any convexity, cf.~\cite{vel-1993}, the one defined here is the rank w.r.t.\ the geodesic convexity. 

If the open neighborhood $N(x)$ of a vertex $x\in V(G)$ induces a complete graph, then  is $x$ is called {\em simplicial}. In other words, $x$ is simplicial if and only if $N(x)$ is convex in $G$. The subgraph of $G$ induced by $S\subseteq V(G)$ is denoted by $G[S]$. The order of a largest complete subgraph of a graph $G$ is $\omega(G)$ and the order of its largest independent set is $\alpha(G)$. A vertex of degree $n(G)-1$ is {\em universal}. The complement of $G$ will be denoted with $\overline{G}$. The {\em join} $G+H$ of graphs $G$ and $H$ is the graph obtained by first taking the the disjoint union of $G$ and $H$, and then   adding all possible edges between vertices of $G$ and $H$. Set finally$[n] = \{1,\ldots,n\}$, where $n\in \mathbb{N}$.

\section{The characterization}    
\label{sec:characterization}

In this section we characterize general position sets in graphs. For this sake the following concepts are needed. 

Let $G$ be a connected graph, $S\subseteq V(G)$, and ${\cal P} = \{S_1, \ldots, S_p\}$ a partition of $S$. Then ${\cal P}$ is \emph{distance-constant} if for any $i,j\in [p]$, $i\ne j$, the value $d(u,v)$, where $u\in S_i$ and $v\in S_j$ is independent of the selection of $u$ and $v$. (We note that in~\cite[p.~331]{kante-2017} the distance-constant partition is called ``distance-regular'', but we decided to rather avoid this naming because distance-regular graphs form a well-established term, cf.~\cite{brouwer-1989}.) If ${\cal P}$ is a distance-constant partition, and $i,j\in [p]$, $i\ne j$, then let $d(S_i, S_j)$ be the distance between the sets $S_i$ and $S_j$, that is, the distance between two arbitrary vertices pairwise from them. Finally, we say that a distance-constant partition ${\cal P}$ is {\em in-transitive} if $d(S_i, S_k) \ne d(S_i, S_j) + d(S_j,S_k)$ holds for arbitrary pairwise different $i,j,k\in [p]$. With these concepts in hand we can characterize general position sets as follows. 

\begin{theorem}
\label{thm:gpsets}
Let $G$ be a connected graph. Then $S\subseteq V(G)$ is a general position set if and only if the components of $G[S]$ are complete subgraphs, the vertices of which form an in-transitive, distance-constant partition of $S$. 
\end{theorem}

\proof
Let $S\subseteq V(G)$. Let $G_1, \ldots, G_p$ be the components of $G[S]$ and let ${\cal P}$ be the partition of $S$ induced by the vertex sets of the components, that is, ${\cal P} = \{V(G_1), \ldots, V(G_p)\}$. To simplify the notation let $V_i = V(G_i)$ for $i\in [p]$, so that ${\cal P} = \{V_1, \ldots, V_p\}$.

Suppose first that $G_1, \ldots, G_p$ are complete subgraphs of $G$ and that ${\cal P}$ forms an in-transitive, distance-constant partition of $S$. We claim that $S$ is a general position set. Assume by the way of contradiction that $S$ contains three vertices $u,v,w$, such that $v$ lies in $I(u,w)$. Since $G_1, \ldots, G_p$ are complete subgraphs, $u$ and $w$ lie in different parts of ${\cal P}$, say $u\in V_i$ and $w\in V_j$, where $i,j\in [p]$, $i\ne j$. Since ${\cal P}$ is distance-constant, we infer that $v\notin V_i$ as well as $v\notin V_j$. Therefore, $v\in V_k$ for some $k\in [p]$, $k\ne i,j$. But then $d(V_i, V_j) = d(V_i,V_k) + d(V_k, V_j)$, a contradiction with the assumption that ${\cal P}$ is an in-transitive partition. 

Conversely, let $S$ be a general position set. If $G_i$ is not complete for some $i\in [p]$, then $G_i$ contains an induced $P_3$, say $uvw$. But this means that $S$ is not a general position set. Hence $G_i$ is   complete for every $i\in [p]$. Next, let $u,v\in V_i$ and $w\in V_j$ for $i,j\in [p]$, $i\ne j$. Since $G_i$ and $G_j$ are complete, $uv\in E(G)$ and hence $|d(u,w) - d(v,w)|\le 1$. Moreover, neither $v$ can be on a $u,w$-geodesic, nor $u$ lies on a $v,w$-geodesic and consequently $d(u,w) = d(v,w)$. Since $u,v$ are arbitrary vertices of $G_i$ and $w$ an arbitrary vertex of $G_j$, this means that $d(V_i, V_j) = d(u,w) = d(v,w)$ is well defined. Consequently, ${\cal P}$ is a distance-constant partition. Finally, ${\cal P}$ must also be an in-transitive partition. If this would not be the case, then there would exist components (complete subgraphs) $G_i$, $G_j$, and $G_k$ of ${\cal P}$ such that $d(V_i, V_k) = d(V_i, V_j) + d(V_j,V_k)$. But this means that if $u\in V_i$, $v\in V_j$, and $w\in V_k$, then $v$ would lie on a $u,w$-geodesics. This contradiction implies that ${\cal P}$ is indeed an in-transitive partition. 
\qed

Theorem~\ref{thm:gpsets} in particular implies some earlier results. First, it immediately implies~\cite[Lemma~3.5]{manuel-2018a} asserting that the set of simplicial vertices of a given graph lies in a general position. Also, setting $d(e,f)$ = $\min\{d(u,x), d(u,y), d(v,x), d(v,y)\}$ for edges $e=uv$ and $f=xy$ of a graph $G$, we obtain: 

\begin{corollary} {\rm \cite[Proposition 4.4]{manuel-2018a}}
Let $G$ be a graph with $\diam(G) \ge 2$. If $F\subseteq E(G)$ is such that $d(e,e') = {\rm diam}(G)$, $e,e'\in F$, $e\ne e'$, then $\gp(G)\ge 2|F|$.  
\end{corollary}

\proof
For $e\in F$ let $x_e$ and $y_e$ be the end-vertices of $e$. Then, having in mind that $\diam(G) \ge 2$, it is straightforward to see that $\{\{x_e, y_e\}:\ e\in F\}$ forms an in-transitive, distance-constant partition.   
\qed

\section{Graphs of diameter $2$}    
\label{sec:diameter-2}

Graphs of diameter $2$ form one of the most interesting classes of graph theory, after all, as it is well-known, almost all graphs have diameter $2$. They are still extensively investigated, the papers~\cite{araujo-2017, bickle-2019, xu-2018} are examples of recent developments on these graphs. In this section we are going to use Theorem~\ref{thm:gpsets} in the case of graphs of diameter $2$. For this sake we denote with $\eta(G)$ the maximum order of an induced complete multipartite subgraph of $\overline{G}$. Note that $K_n$ is complete multipartite, and that $\eta(K_n) =1$ and $\omega(K_n)=n$. 

\begin{theorem}
\label{thm:diameter2}
If $\diam(G) = 2$, then $\gp(G) = \max\{\omega(G), \eta(G)\}$. 
\end{theorem}

\proof
Since the vertices of an arbitrary complete subgraph of a graph $G$ form a general position set of $G$, we have $\gp(G) \ge \omega(G)$. Suppose $H$ is a complete multipartite subgraph of $\overline{G}$. Then in $G$ the vertices of $H$ induce a disjoint union of complete graphs. Since  $\diam(G) = 2$, the vertices of these complete subgraphs clearly form an in-transitive, distance-constant partition. Hence by Theorem~\ref{thm:gpsets}, $\gp(G) \ge \eta(G)$. Therefore, $\gp(G) \ge \max\{\omega(G), \eta(G)\}$. 

Let now $S$ be a set of vertices in a general position in $G$. Then by Theorem~\ref{thm:gpsets} the components of $G[S]$ are complete subgraphs, the vertices of which form an in-transitive, distance-constant partition of $S$. If there is only one such component, then $|S|\le \omega(G)$, and if there are at least two components, then $|S|\le \eta(G)$. Hence, $\gp(G) \le \max\{\omega(G), \eta(G)\}$. 
\qed

If $P$ is the Petersen graph, then $\omega(P) = 2$ and $\eta(P) = 6$, hence by Theorem~\ref{thm:diameter2} we have $\gp(P) = 6 = \eta(P)$. Let further $G_{n,k}$, be the graph obtained from $K_n$ and one more vertex that is adjacent to $k+1$ vertices of $K_n$, where $n\ge 3$ and $1\le k+1 < n$. Then $\omega(G) = n$ and $\eta(G) = n-k$, so that $\gp(G_{n,k}) = n = \omega(G)$. These examples show that the values from the maximum in Theorem~\ref{thm:diameter2} are independent. 

Cographs form an important class of graphs that is still extensively investigated,~\cite{allem-2018, tsujie-2018} is a selected couple of recent studies. Recall that $G$ is a {\em cograph} if $G$ contains no path $P_4$ as an induced subgraph. These graphs were independently introduced several times and can be characterized in many different ways, see~\cite{corneil-1981}. In particular, a graph is a cograph iff it can be obtained from $K_1$ by means of the disjoint union and join of graphs. Note that this implies that every connected cograph of order at least $2$ is the join of at least two smaller connected cographs. This implies that $\diam(G) \le 2$ holds for any connected cograph $G$. 

\begin{theorem}
\label{thm:cographs}
If $G$ be a connected cograph, then $\gp(G) = \max\{\omega(G), \alpha(G)\}$.
\end{theorem}

\proof
If $G = K_n$, then $\gp(K_n) = n = \max\{\omega(K_n), \alpha(K_n)\}$. Hence assume in the rest that $G$ is a connected cograph with $\diam(G) = 2$.

We claim that $\eta(G) = \alpha(G)$ and proceed by induction on the order of $G$. The assertion is clear if $n(G) = 3$ (in other words, for $G=P_3$). Assume now that $G$ is a connected cograph with $\diam(G) = 2$ and $n(G)\ge 4$. Then $G = G_1 + \cdots +G_k$, where $k\ge 2$ and $G_i$, $i\in [k]$, are connected cographs. Since for arbitrary graphs $X$ and $Y$ we have $\alpha(X+Y) = \max\{\alpha(X), \alpha(Y)\}$ and $\eta(X + Y) = \max\{\eta(X), \eta(Y)\}$, we get, by the induction assumption, that 
\begin{eqnarray*}
\label{eq:alpha-eta}
\eta(G) = \max\{\eta(G_1),\ldots, \eta(G_k)\} = \max\{\alpha(G_1),\ldots, \alpha(G_k)\} = \alpha(G)\,. 
\end{eqnarray*}
The result now follows from Theorem~\ref{thm:diameter2}. 
\qed

If $G$ is a cograph, then $\alpha(G)$ and $\omega(G)$ can be determined in polynomial time, cf.~\cite{corneil-1981, orlovich-2011}. Hence Theorem~\ref{thm:cographs} implies that the general position problem is polynomial on connected cographs. Since the general position function of a graph is clearly additive on its components, the general position problem is thus polynomial on all cographs.

\medskip
Suppose that $G$ is not complete and that it contains at least one universal vertex. Then $\diam(G) = 2$ and Theorem~\ref{thm:diameter2} applies. For this situation we get: 

\begin{corollary}
\label{cor:Delta=n-1}
Let $G$ be a non-complete graph, $U\ne \emptyset$ the set of its universal vertices, and let $U' = V(G) - U$. Then 
$$\gp(G) = \max\{|U| + \omega(G[U']), \eta(G[U'])\}\,.$$ 
\end{corollary}

\proof
Since $U$ contains universal vertices, every largest complete subgraph of $G$ contains $U$. Hence $\omega (G) = |U| + \omega(G[U'])$. In $\overline{G}$ every vertex from $U$ is isolated. Hence every induced complete multipartite subgraph of $\overline{G}$ with at least two parts contains only vertices from $U'$. Thus $\eta(G) =  \eta(G[U'])$. 
\qed

Every graph $G$ can be represented as the graph obtained from $K_{n(G)}$ by removing appropriate edges. To present examples how Corollary~\ref{cor:Delta=n-1} can be applied, let us use the notation $K_n - E(H)$ for the graph obtained from $K_n$ in which we consider $H$ as its subgraph, and then deleting the edges of $H$ from $K_n$. Then we have the following formulas that can be easily deduced from Corollary~\ref{cor:Delta=n-1}, where $W_k$ denotes the wheel graph of order $k$, that is, the graph obtained from $C_{k-1}$ by adding an additional vertex and making it adjacent to all the vertices of $C_{k-1}$. 

\begin{itemize}
\item $\gp(K_n - E(K_k)) = \max \{ k, n- k + 1 \}$, where $2\leq k < n$.
\item $\gp(K_n - E(K_{1,k})) = \max \{ k+1, n-1 \}$, where $ 2\leq k < n$. 
\item $\gp(K_n - E(P_k)) = \max \{3, n- k + \lceil\frac{k}{2}\rceil \}$, where $3\leq k< n$. 
\item $\gp(K_n - E(K_{r,s})) = \max \{ r+s, n-r \}$, where $ 2\leq r\leq s$ and $r +s < n$. 
\item $\gp(K_n - E(W_k)) = \max \{3, n- k + \lfloor\frac{k-1}{2}\rfloor \}$, where $ 5\leq k< n$. 
\item $\gp(K_n - E(C_k)) =
   \left\{
     \begin{array}{ll}
       \max \{ 3, n- k + \lfloor\frac{k}{2}\rfloor\}, & 5 \le k < n; \\
       \max \{ 4, n-2 \}, & k = 4.
     \end{array}
   \right.$
\end{itemize}

\section{Bipartite graphs and their complements}    
\label{sec:bipartite}

For bipartite graphs we have the following result. 

\begin{theorem} 
\label{thm:bipartite}
If $G$ is a connected, bipartite graph on at least $3$ vertices, then $\gp(G) \leq \alpha(G)$.  Moreover, if $\diam(G) \in \{2,3\}$, then $ \gp(G)=\alpha (G)$.
\end{theorem}

\proof
If $G$ is a connected graph, then $\gp(G) = 2$ if and only if $G = P_n $ $(n\ge 2)$ or $G = C_4$, see~\cite[Theorem 2.10]{ullas-2016}. Thus $\gp(G) \leq \alpha (G)$ holds in these cases. (In the case of $P_n$ the inequality  holds because we have assumed that $n\ge 3$ and so $\alpha(P_n) \ge 2$.) In the rest we may thus assume that $G$ is a connected, bipartite graph that is neither a path nor $C_4$, so $\gp(G) \geq 3$. 

Let $S$ be a gp-set of $G$ and let $ S_1, \ldots, S_k$ be the components of $G[S]$. As $\gp(G) \geq 3$ and $G$ is bipartite, Theorem~\ref{thm:gpsets} implies that $k\ge 2$. Also, since $G$ is bipartite, $|S_i|\in [2]$ for $i\in [k]$. We claim that actually $|S_i| = 1 $ for every $i\in [k]$. Suppose on the contrary that, w.l.o.g., $S_1 = \{u, v\}$. Let $w\in S_2$. Then, since $uv\in E(G)$ and $G$ is bipartite, $|d(u,w) - d(v,w)| = 1$, which means that either $v$ lies on a $u,w$-geodesic or $u$ lies on a $v,w$-geodesic. This contradiction proves the claim, that is, $S$ is an independent set. We conclude that $\gp(G) \leq \alpha(G)$.

Assume now that $\diam(G) = 2$. The complement of an independent set $I$ of a graph $G$ induces the complete graph $K_{|I|}$ which is an instance of a complete multipartite graph. Hence $\alpha(G) \le \eta(G)$ and Theorem~\ref{thm:diameter2} implies that $\alpha(G)\le \gp(G)$ holds for a graph $G$ of diameter $2$. 

Assume finally that $\diam(G) = 3$. Then we recall from~\cite[Corollary 4.3]{manuel-2018a} that every independent set is a general position set. Hence $\alpha(G)\le \gp(G)$ holds also in this case. 
\qed

If $G$ is bipartite, $\gp(G)$ can be arbitrary smaller than $\alpha(G)$. Consider first the paths $P_n$, $n\ge 2$, for which we have $\gp(P_n) = 2$ and $\alpha(P_n) = \lceil n/2\rceil$. We also note that none of the two assertions of Theorem~\ref{thm:bipartite} need hold if $G$ is not bipartite. To see this, consider again the Petersen graph $P$. (Of course, $\diam(P) = 2$.) As already noticed, $\gp(P) = 6$, while $\alpha(P) = 4$. For a corresponding example of diameter $3$ just add a pendant vertex to $P$. 

We now turn our attention to complements of bipartite graphs for which some preparation is needed. If $G = (V(G), E(G))$ is a bipartite graph and $V(G) = A\cup B$ is its bipartition, then we will write $G$ as triple $(A, B, E(G))$.  If $G = (A, B, E(G))$ is a bipartite graph, then let $M_G$ be the set of vertices of largest possible degree, more precisely,  
$$M_G = \{u\in A:\ {\rm deg}(u) = |B|\}\, \cup\, \{u\in B:\ {\rm deg}(u) = |A|\}\,.$$
Let $\psi(G)$ be the maximum order of an induced complete bipartite subgraph of $G$. Note that if $G$ is a bipartite, but not complete bipartite, then $\diam(\overline{G})\le 3$. Now we can formulate: 

\begin{theorem}
\label{thm:bip-complement}
If $G = (A, B; E(G))$ is a bipartite graph, then 
$$\gp(\overline{G}) =
   \left\{
     \begin{array}{ll}
       n(G), & \hspace*{-1.1cm}\diam(\overline{G}) \in \{1, \infty\}; \\ 
       \max \{\alpha(G), \psi(G)\}, & \hspace*{-0.1cm} \diam (\overline{G}) = 2; \\
\max \{\alpha(G), \psi(G\setminus (M_G\cap A)), \psi(G \setminus (M_G\cap B)), |M_G|\}, & \hspace*{-0.1cm} \diam(\overline{G}) = 3.
\end{array}
\right.$$
\end{theorem}

\proof
Let $G = (A, B; E(G)) $ be a bipartite graph. Then $\overline{G}$ is disconnected iff $G$ is complete bipartite. In this case we have $\diam(\overline{G}) = \infty$ and $\overline{G}$ is a disjoint union of $K_{|A|}$ and $K_{|B|}$. Therefore, $\gp(\overline{G}) = |A| + |B| = n(G)$. Further, $\diam(\overline{G}) = 1$ iff $G$ is edge-less, hence again $\gp(\overline{G}) = n(G)$. If $\diam (\overline{G}) = 2$, then by Theorem~\ref{thm:diameter2} we have $\gp(\overline{G}) = \max \{\omega(\overline{G}), \eta(\overline{G})\}$. Since $\omega(\overline{G}) = \alpha(G)$ and $\eta(\overline{G}) = \psi(G)$, the assertion for the diameter $2$ follows. 

In the rest we may thus assume $|A|\ge 2$, $|B|\ge 2$, and $\diam(\overline{G}) = 3$. Note that if $u\in M_G\cap A$, then $u$ has no neighbor in $B$ and if $u\in M_G\cap B$, then $u$ has no neighbor in $A$. Consequently, in $\overline{G}$ two vertices are at distance $3$ if and only if one lies in $M_G\cap A$ and the other in $M_G\cap B$. Since $\diam(\overline{G}) = 3$ it follows that $M_G\cap A \ne \emptyset$ and $M_G\cap B \ne \emptyset$. 

Consider a set $T$ in general position in $\overline{G}$ and set $T_A = T\cap A$, $T_B = T\cap B$. If $T$ has at least one vertex in $M_G\cap A$, say $x$, and at least one vertex in $M_G\cap B$, say $y$, then every vertex from $(A\cup B)\setminus M_G$ lies on a $x,y$-geodesic. Therefore, $T\subseteq M_G$. This means that $|T|\le |M_G|$. Suppose next that $T\cap (M_G\cap A) = \emptyset$. If there is an edge between a vertex from $T_A$ and a vertex from $T_B$, then $T$ must induce a clique and hence $|T|\le \omega(\overline{G}) = \alpha(G)$. Otherwise, in view of Theorem~\ref{thm:gpsets}, the vertices from $T_A$ and from $T_B$ are pairwise at distance $2$. But then $T$ induces a complete bipartite graph in $G\setminus (M_G\cap A)$ and therefore $|T|\le \psi(G\setminus (M_G\cap A))$. Analogously, if $T\cap (M_G\cap B) = \emptyset$ then we get that  $|T|\le \alpha(G)$ or $|T|\le \psi(G\setminus (M_G\cap B))$. In summary, 
$$\gp(\overline{G}) \le \max \{\alpha(G), \psi(G\setminus (M_G\cap A)), \psi(G \setminus (M_G\cap B)), |M_G|\}\,.$$
On the other hand, we clearly have $\gp(\overline{G}) \ge \omega(\overline{G}) = \alpha(G)$. Note next that each vertex from $M_G$ is simplicial in $\overline{G}$ and consequently $\gp(\overline{G}) \ge |M_G|$. Finally, an induced complete bipartite graph in $G\setminus (M_G\cap A)$ as well as in $G\setminus (M_G\cap B)$ corresponds to a disjoint union of cliques in $\overline{G}$ which form an in-transitive, distance constant partition (with constant $2$). Hence we also have $\gp(\overline{G}) \ge \psi(G\setminus (M_G\cap A))$ and $\gp(\overline{G}) \ge \psi(G \setminus (M_G\cap B))$. 
\qed 

If $n\ge 5$, then $\diam(\overline{P_n}) = 2$ and for $n\ge 7$ we have $\psi(P_n) = 3 <  \lceil n/2\rceil = \alpha(P_n)$.  Let next $G_n$ be a bipartite graph with the bipartition $A = \{x_1,\ldots, x_n, a_1, a_2\}$ and $B = \{y_1,\ldots, y_n, b_1, b_2\}$, where vertices $(A\cup B)\setminus \{a_1,a_2,b_1,b_2\}$ induce a complete bipartite graph $K_{n,n}$, and the remaining edges of $G_n$ are $a_1y_1$, $a_2y_2$, $b_1x_1$, and $b_2x_2$. For $n\ge 3$ we have $\psi(G_n) = 2n > n+2 = \alpha(G_n)$. As $\diam(\overline{G_n}) = 2$, these two examples demonstrate that the values in Theorem~\ref{thm:bip-complement} are independent in the case $\diam(\overline{G}) = 2$. 

Let $H(n,m,s,t)$, $n,m,s,t\ge 2$, be a bipartite graph with the bipartition $A=A_1\cup A_2\cup A_3$ and $B = B_1\cup B_2\cup B_3$, where $|A_1| = n$, $|B_1| = m$, $|A_2| = |B_3| = s$, and $|A_3| = |B_2| = t$. The vertices in $(A_1\cup A_2)\cup (B_1\cup B_2)$ induce $K_{n+s,m+t}$, the vertices in $A_2\cup B_3$ induce $K_{s,s}$ and the vertices in $A_3\cup B_2$ induce $K_{t,t}$. These are all the edges of $H(n,m,s,t)$. Assume that $n\le m$ and set $H = H(n,m,s,t)$. Then $M_H \cap A = A_2$ and $M_H \cap B = B_2$ and we have: 
\begin{itemize}
\item $|M_H| = s + t$,
\item $\alpha(H) = m + s + t$,
\item $\psi(H \setminus A_2) = \max\{m + n + t, 2t\}$, and 
\item $\psi(H \setminus B_2) = \max\{m + n + s, 2s\}$.  
\end{itemize}
It is now clear that the parameters $n$, $s$, and $t$ can be selected such that exactly one of $\alpha(H)$, $\psi(H \setminus A_2)$, and $\psi(H \setminus B_2)$ is strictly larger than the other two (as well as bigger than $|M_H|$). Note finally that $\diam(\overline{H}) = 3$. 

To see that $|M_G|$ can be strictly larger than the other three terms from Theorem~\ref{thm:bip-complement} when $\diam(\overline{G}) = 3$, consider the edge deleted complete bipartite graph $K = K_{n,n}-e$. Then $\diam(\overline{K}) = 3$, and $|M_K| = 2n-2$.   

In the rest we present the general position number for some natural families of bipartite complements. 

In~\cite[Theorem 2.5]{ullas-2016} and in~\cite[Corollary 3.7]{manuel-2018a} it was independently observed that the gp-number of a tree $T$ is the number of its leaves. (Actually, the set of leaves is the unique gp-set of $T$.) For the complements of trees we have: 

\begin{corollary}
It $T$ is a tree, then 
$\gp(\overline{T}) = \max \{\alpha(T), \Delta(T) +1\}$.
\end{corollary}

\proof
Let $T = (A,B; E(T))$. 

If $\diam(T)\leq 2$, then it is clear that $T$ is a star. Hence $\diam(\overline{T})=\infty$. By  Theorem~\ref{thm:bip-complement} we thus have $\gp(\overline{T}) = n(T)= \Delta(T) +1 $. 

If $\diam(T)=3$, then it is straightforward to see that $T$ is isomorphic to a double star. Therefore, $| M_T\cap A | = 1$ and $| M_T \cap B| = 1$. Thus $|M_T| = 2$, $\alpha(T) = |n(T)|- 2$, and $\psi(T\setminus (M_T\cap A))= |A|;$ $ \psi(T\setminus (M_T\cap B)) =|B|$. Since $|n(T)|\ge 4$ and $|A|\ge 2$ and $|B|\ge 2$, we have $\gp(\overline{T}) = |n(T)|-2 = \alpha(T)$. 

Let finally $\diam (T) \geq 4$. Then from~\cite[Lemma 2.2]{senbagamalar-2014} we deduce that $\diam(\overline{T}) = 2$. Since $T$ has no cycles, we have that $\psi(T)$ is the order of a maximum induced star, that is, $\psi(T)=\Delta(T) +1$. Thus $\gp(T) = \max \{\alpha(T), \Delta(T) +1\}$.
\qed

The {\em Cartesian product} $G\cp H$ of graphs $G$ and $H$ is defined as follows. $V(G\cp H) = V(G) \times V(H)$. As for the edges, $(g,h)$ and $(g',h')$ are adjacent if (i) $g=g'$ and $hh'\in E(H)$, or (ii) $h=h'$ and $gg'\in E(G)$. (See the book~\cite{hik-2011}.) In~\cite[Theorem 3.1]{manuel-2018b} it was proved that $\gp(P_\infty \cp P_\infty) = 4$, where $P_\infty$ is the two-way infinite path. If follows from this result that if $n,m\ge 3$, then $\gp(P_n\cp P_m) = 4$ as well. For the complements of these grids we have: 

\begin{corollary}
\label{cor:grids}
If $n,m\ge 2$, then 
$$ gp(\overline{P_n\cp P_m}) =
   \left\{
     \begin{array}{ll}
       4, &  n = m = 2;\\
       \left\lceil \frac{n}{2}\right\rceil \left\lceil \frac{m}{2}\right\rceil + \left\lfloor \frac{n}{2}\right\rfloor \left\lfloor \frac{m}{2}\right\rfloor, & {\rm otherwise}\,.
     \end{array}
   \right.$$
\end{corollary}

\proof
$\overline{P_2\cp P_2} = \overline{C_4}$, hence the assertion holds for $n = m = 2$. $\overline{P_2\cp P_3}$ is the graph obtained from the $6$-cycle $v_1v_2\ldots v_6$ by inserting into it the edges $v_1v_3$ and $v_4v_6$. The assertion then follows immediately. 

Suppose in the rest that $n,m \ge 3$ and set $G = P_n\cp P_m$. Applying~\cite[Lemma 2.2]{senbagamalar-2014} once more we get that $\diam(\overline{G}) = 2$. Hence by  Theorem~\ref{thm:bip-complement} we see that $\gp(G)=\max \{\alpha(G), \psi(G)\}$. Since the only induced complete bipartite subgraphs in $P_n\cp P_m$ are isomorphic to $K_{2,2}$ or $K_{1,r}$, $r\in [4]$, we get $\gp(G) = \alpha(G)$. The conclusion of the theorem now follows because $\alpha(P_n\cp P_m) = \left\lceil \frac{n}{2}\right\rceil \left\lceil \frac{m}{2}\right\rceil + \left\lfloor \frac{n}{2}\right\rfloor \left\lfloor \frac{m}{2}\right\rfloor$, a result that can be deduced from~\cite[Theorem 4.2]{klavzar-2005}. 
\qed

Using parallel arguments as in the proof of Corollary~\ref{cor:grids} we also get the general position number of the complements of hypercubes. 

\begin{corollary} 
If $k\ge 3$, then $gp(\overline{Q_k}) = 2^{k-1}$. 
\end{corollary}

\section*{Acknowledgements}

S.K.\ acknowledges the financial support from the Slovenian Research Agency (research core funding P1-0297 and projects J1-9109, N1-0095). E.J.T.\ acknowledges the financial support from the University of Kerala for providing University JRF.


\end{document}